\documentclass[11pt]{amsart}
\usepackage[dvipsnames]{xcolor}
\usepackage{amssymb,verbatim}
\usepackage{amsmath,amsfonts,enumitem,bm}
\usepackage[mathscr]{euscript} % make a nice \Net
\usepackage{amsthm}
\usepackage{url}
\usepackage{graphicx} %to include figure
\usepackage{enumitem} %to use nice enumeration with changed indent

\usepackage{hyperref} % to make links for lemmas and theorems

\hypersetup{colorlinks} %so the links look nicer
\usepackage{bbm} % for the characteristic function 1
\usepackage{bm} % for bold greek letter
\usepackage{mathrsfs} %to use \mathscr
\usepackage{cancel} %to make sout to cross things out
\usepackage{ytableau} %to draw tableaux
\usepackage{mathtools}
\usepackage{faktor} %%nice quotients

\usepackage[colorinlistoftodos]{todonotes}

\RequirePackage{cleveref}
\usepackage{hypcap}
\hypersetup{colorlinks=true, citecolor=darkblue, linkcolor=darkblue}
\definecolor{darkblue}{rgb}{0.0,0,0.7}

\definecolor{darkred}{rgb}{0.68,0,0}

\definecolor{darkgreen}{rgb}{0,.38,0}

\setlist[enumerate]{
	label=\textnormal{({\roman*})},
	ref={\roman*}}

\makeatletter
\def\th@plain{%
	\thm@notefont{}% same as heading font
	\itshape % body font
}
\def\th@definition{%
	\thm@notefont{}% same as heading font
	\normalfont % body font
}
\makeatother

\newtheorem{thm}{Theorem}[section]
\newtheorem{lemma}[thm]{Lemma}

\newtheorem*{claim*}{Claim}

\newtheorem{conv}[thm]{Convention}

\theoremstyle{definition}

\newtheorem{rem}[thm]{Remark}

\numberwithin{figure}{section}
\numberwithin{equation}{section}

\def\bu{\bullet}

\def\zz{\mathbb Z}
\def\nn{\mathbb N}
\def\cc{\mathbb C}

\def\ov{\overline}

\def\Ga{\Gamma}

\def\la{\lambda}
\def\ga{\gamma}

\def\ep{\ve}
\def\al{\alpha}

\def\ve{\varepsilon}

\def\vk{\varkappa}

\def\cT{\mathcal T}

\def\ssu{\subset}

\def\<{\langle}
\def\>{\rangle}

\def\0{{\mathbf 0}}

\def\TT{{\mathbb{T}}}

\def\.{\hskip.06cm}
\def\ts{\hskip.03cm}

\def\lra{\leftrightarrow}

\def\pt{\partial}

\def\bx{{\textbf{\textit{x}}}}

\def\bw{\textbf{\textit{w}}}

\def\bb{\textbf{\textbf{b}}}

\def\bal{{\boldsymbol{\alpha}}}

\def\bb{{\boldsymbol{b}}}

\def\di{{\small{\ts\diamond\ts}}}

\def\ze{{\zeta}}

\newcommand{\supp}{\mathrm{NF}}

\def\.{\hskip.06cm}
\def\ts{\hskip.03cm}

\def\nin{\noindent}

\def\bal{\mathbf{\al}}

\newcommand{\textsu}[1]{\textup{\textsf{#1}}}
\usepackage{indentfirst}
\DeclareTextSymbolDefault{\ae}{T1}

\newcommand{\ComCla}[1]{\textup{\textsu{#1}}}
\newcommand{\sharpP}{\ComCla{\#P}}
\newcommand{\SP}{\ComCla{\#P}}

\newcommand{\GapP}{\ComCla{GapP}}

\renewcommand{\P}{\ComCla{P}}

\newcommand{\AM}{\ComCla{AM}}
\newcommand{\coAM}{\ComCla{coAM}}

\def\SP{\sharpP}

\def\GRH{\textup{\sc GRH}}

\newcommand{\inv}{\operatorname{{\rm inv}}}
\newcommand{\Des}{\operatorname{{\rm Des}}}

\newcommand{\tinysquare}{\scriptscriptstyle{\square}}

\newcommand{\Sch}{\mathfrak{S}} %% Macro for Schubert polynomials

\DeclareMathOperator{\one}{\mathbf{1}} % vector with value 1 at the first entry 0 everything else
\newcommand{\onMark}{\ast}
\newcommand{\offMark}{\circ}

\begin{document}

\title[Signed puzzles for Schubert coefficients]{Signed puzzles for Schubert coefficients}

\author[Igor Pak \. \and \. Colleen Robichaux]{Igor Pak$^\star$  \. \and \.  Colleen Robichaux$^\diamond$}

\thanks{\today}
\thanks{\thinspace ${\hspace{-.45ex}}^\star$Department of Mathematics,
UCLA, Los Angeles, CA 90095, USA. Email:  \texttt{pak@math.ucla.edu}}
\thanks{\thinspace ${\hspace{-.45ex}}^\diamond$Department of Mathematics,
UC Davis, Davis, CA 95616, USA. Email:  \texttt{robichaux@ucdavis.edu}}

\begin{abstract}
We give a signed puzzle rule to compute Schubert coefficients.  The rule is
based on a careful analysis of Knutson's recurrence \cite{Knu03}.  We use the
rule to prove polynomiality of the sums of Schubert coefficients with bounded
number of inversions.
\end{abstract}

%%%%%%%%%%%%%%%%%%%%%%%%%%%%%%%%%%%%%%%%%%%%%%%%%%%%%%%%%%%%
%  TITLE PAGE information
%%%%%%%%%%%%%%%%%%%%%%%%%%%%%%%%%%%%%%%%%%%%%%%%%%%%%%%%%%%%

\maketitle

\section{Introduction}\label{s:intro}

% \subsection{Main result} \label{s:intro-main}
%
Schubert coefficients are extremely well studied yet deeply mysterious
numbers which play a central role in Schubert calculus.  A major open
problem asks for a combinatorial interpretation of the coefficients
\cite[Problem~11]{Sta00}.  Puzzles are finite tilings with edge-labeled
equilateral triangles which enumerate the desired numbers.  Such puzzle
rules were discovered for many special cases and for closely related problems;
we refer to \cite{Knu22} for an extensive overview.

While a manifestly positive combinatorial interpretation remains elusive,
signed combinatorial interpretations are also of interest for various
applications, see a discussion in \cite{Pak-OPAC}. In \cite{PR24a},
the authors present signed combinatorial interpretations for a
wide range of structure constants in algebraic combinatorics,
including Schubert coefficients. For Schubert coefficients
and their generalizations, several such formulas are known in the
literature, see a discussion in~$\S$\ref{ss:finrem-other}.
Unfortunately, neither of these signed combinatorial interpretations
can be extended to a signed puzzle rule.

In this paper we present a signed puzzle rule for Schubert coefficients.
Similar (signed) puzzle rules already exist in a few special cases,
see \cite{KZ21,KZ23}.
Our result is the first signed puzzle rule which holds in full generality.

\smallskip

\begin{thm}\label{t:main}
For every integer $n$, let \ts $\cT_n$ \ts be a set of \ts $O(n^9)$ \ts
puzzle pieces defined in Section~\ref{s:puzzle}.
Let \ts $u,v,w\in S_n$ \ts be permutations with \ts $\inv(u)+\inv(v)=\inv(w)$,
and denote \ts $\ell=\binom{n}{2}-\inv(u)$.
Let  \ts $\Ga=\Ga(u,v,w)$ \ts be an \ts $n \times \ell$ \ts
parallelogram region with indicators and labels defined in Section~\ref{s:puzzle}.
Then the number of signed puzzles of \ts $\Ga$ \ts with~$\ts \cT_n$ \ts
is the Schubert coefficient \ts $c_{u,v}^w\ts$.
\end{thm}

\smallskip

The starting point of our construction is a special case of \emph{Knutson's
recurrence} given in \cite{Knu03}, see Section~\ref{s:knutson} below.
% Curiously,
Knutson's recurrence is an advanced extension of an earlier
paper \cite{Knu01}.
Note that we do not use the equivariant variables and consider the
results only in type~$A$.
See also Yong's implementation of Knutson's recurrence \cite{Yong06}.

\smallskip

The proof of Theorem~\ref{t:main} is completely combinatorial,
and uses only basic notions from Schubert calculus.
The construction of the puzzles is given in Section~\ref{s:puzzle};
it is somewhat involved but richly illustrated.
The proof of the theorem is then given in Section~\ref{s:proof},
and a large example is given in Section~\ref{s:ex}.
We also give the following unusual application of the signed puzzle rule:

\smallskip

\begin{thm}  \label{t:app}  Fix~$k$, and let
$$
\ga_k(n) \, := \, \sum_{u,v,w\ts\in\ts S_n \. : \. \inv(w)=k} \. c^w_{u,v}\..
$$
Then \ts $\ga_k$ \ts is a polynomial in~$n$.
\end{thm}

\smallskip

Since we have \ts $u,v\leqslant w$ \ts in Bruhat order for the nonzero terms
in the summation, this gives
\ts $\inv(u),\inv(v)\le k$.  Thus the total number of triples \ts $(u,v,w)$ \ts
in the summation above is \emph{at most} \ts $\binom{n}{2}^{3k} \le n^{6k}$.
To bound the Schubert coefficients in the summation, note that
 \ts $|\supp(u)\cup\supp(v)|\le 4k$, where \ts $\supp(w):=\{i\in[n] \, : \, w(i)\neq i\}$
 \ts denotes the set of non-fixed points in~$w$.   Stanley's upper bound
 in \cite[$\S$5]{Sta17} gives  \ts $c^w_{v,w} \le 2^{(4k)^2}$.
 Therefore, \ts $\ga_k(n)= O_k(n^{6k})$.  However, a priori there is no reason
 to believe that the sum \ts $\ga_k$ \ts is polynomial in~$n$.  See~$\S$\ref{ss:finrem-inv}
 for further discussion.

The proof of Theorem~\ref{t:app} is given in Section~\ref{s:app}.  It uses
technical details of the puzzle construction in the proof of Theorem~\ref{t:main}
and a geometric argument in Ehrhart theory.  We conclude with final remarks
in Section~\ref{s:finrem}, where we give further comments on the nature of
signed puzzle rules.
% We note that it is straightforward
% to show that \ts $\ga_k$ \ts is a quasi-polynomial; proving polynomiality
% is more difficult and occupies much of the proof.

\medskip

\section{Standard definitions and notation} \label{s:not}
We refer to \cite{Mac91,Man01} for the background on Schubert calculus,
to \cite{Ful97,Sta99} for definitions and standard results in algebraic combinatorics,
and to \cite[$\S$14]{GS87} and \cite{vEB97} for basic results on Wang tilings.
See \cite{BR15} for the introduction to Ehrhart theory of rational
polyhedra, and \cite{Bar97} for a concise survey on the subject.

For two functions \ts $f,g: \nn\to \nn$, we use \ts $f = O(g)$ \ts
if there is a universal constant $C>0$ such that \ts $f(n)\le Cg(n)$ \ts
for all $n\in \nn$.  For two functions \ts $f,g: \nn^2\to \nn$, we
use \ts $f = O_k(g)$ \ts if for all \ts $k \in \nn$ \ts there is
a constant \ts $C(k)>0$ \ts such that \ts $f(n,k)\le C(k) \ts g(n,k)$ \ts
for all \ts $n\in \nn$.

Let \ts $[n]:=\{1,\ldots,n\}$ \ts and \ts $\<n\>:=[n]\cup \{-\}$, where we view ``$-$'' as the \emph{blank element}.
Our notation for permutations will simplify if the context is clear,
e.g.\ we write $4123$ to mean a permutation $(4,1,2,3)$ in~$S_4$\ts.
We think of multiplication on the right as the action on positions, e.g.\
$4123 \cdot 2134 = 1423$.
Let \ts $\inv(w):=\{(i,j) \, : \,i<j,\. w(i)>w(j)\}$ \. denote
the \emph{number of inversions} in~$w$.

For a permutation $w\in S_n$ we say that $i\in [n-1]$ is an \emph{ascent} \ts if
$w(i)<w(i+1)$.  Otherwise, $i$ is a \emph{descent}.   Let \ts $\Des(w)$ \ts
denote the set of descents in~$w$. Let \ts $\bw_\circ=(n,n-1,\ldots,1)$ \ts denotes
the \emph{long permutation}, so \ts $\Des(\bw_\circ)=[n-1]$.  We write $\one$ for the
identity permutation $(1,2,\ldots,n)$.  We use \ts $t_{ij}:=(i,j)$ \ts
to denote a transposition in~$S_n\ts$, and let \ts $s_i:=(i,i+1)$ \ts denote an
adjacent transposition.  By definition, if \ts $i$ \ts a descent in~$w$, then $i$
is an ascent in~$ws_i$ and vice versa.

A \emph{puzzle piece} $\tau$ is a region (tile) in a triangular grid $\TT$ with certain
labels/indicators on the boundary.  In this paper all puzzle pieces will be
unit triangles.  For a collection $\cT$ of puzzle pieces and a
region $\Ga$ in~$\TT$, a \emph{puzzle $T$ \ts of \ts $\Ga$ \ts with~$\ts \cT$} \ts
is a tiling of \ts $\Ga$ \ts with copies of puzzle pieces \ts $\tau\in \cT$
(up to parallel translation), such that the
labels/indicators match along all common edges between the puzzle pieces and
along the boundary of~$\Ga$.  A \emph{signed puzzle} is a puzzle $T$ with a sign function $s(T)\in \{\pm 1\}$.
The \emph{number of signed puzzles} of $\Ga$ with $\cT$ is the sum of signs $s(T)$
over all puzzles $T$ of $\Ga$ with~$\cT$.

\medskip

\section{Schubert coefficients} \label{s:schubert}

Although we will not use Schubert polynomials, for completeness 
we include the original definition due to Lascoux and Sch\"utzenberger \cite{LS82}.
Let
$$
\Sch_{\bw_\circ} \. : = \. x_1^{n-1} x_2^{n-2} \. \cdots  \. x_{n-1}\..
$$
%A permutation $w\in S_n$ is said to have a {descent} at~$i$, if \ts $w(i)> w(i+1)$.
%Denote by \ts $\Des(w)$ \ts the \emph{set of descents} \ts of~$w$, and by \ts $\des(\si):=|\Des(\si)|$
%\ts the \emph{number of descents}. 
Define the \emph{divided difference operator}
$$
\pt_i F \, := \, \frac{F - s_i F }{x_i -x_{i+1}}\,,
$$
where the transposition \ts $s_i$ \ts acts on \ts $F\in \cc[x_1,\ldots,x_n]$ \ts 
by transposing the variables.  Let
$$
\Sch_{ws_i} \. := \. \pt_i \ts \Sch_w \quad \text{for all} \quad i \in \Des(w), 
$$
and define all Schubert polynomials recursively.  It follows
that \ts $\Sch_w\in \zz[\bx]$ \ts are homogeneous polynomials of degree \ts $\inv(w)$.
%Here \. $\inv(w):=\{(i,j) \, : \,i<j,\. w(i)>w(j)\}$ \. is the number of inversions in~$w$.

Finally, \emph{Schubert coefficients} \ts are defined as structure constants:
% {\small
$$
\Sch_u \cdot \Sch_v \, = \, \sum_{w \ts \in \ts S_\infty} \. c^w_{u,v} \. \Sch_w 
\quad \text{for all} \quad u,v\in S_n\ts.
$$

\medskip

\section{Knutson's recurrence} \label{s:knutson}
%
% Let \ts $\inv(w):=\{(i,j) \, : \,i<j,\. w(i)>w(j)\}$ \ts denote the \emph{number of inversions} in~$w$.
It is well known that
\ts $c^w_{u,v}=0$ \ts unless the \emph{dimension equation} \ts holds:
$$
(\oplus) \qquad
\inv(u) \. + \. \inv(v) \ = \ \inv(w)\ts.
$$
Thus we only consider coefficients satisfying $(\oplus)$ in Theorem~\ref{t:main}. Additionally,
$c^w_{u,v}=0$ \ts unless \ts $u\leqslant w$ \ts in Bruhat order.
The following result is a special case of Knutson's recurrence in type~$A$,
adapted in the notation above.

\smallskip

\begin{lemma}[{\rm \emph{Knutson's recurrence} \cite{Knu03}}{}]\label{l:Knutson}
Let \ts $u,v,w\in S_n$ \ts and suppose \ts $i \notin \Des(u)$.  There are four
cases:

\smallskip
{\small $(0)$} If  \ts $i \notin \Des(v)$ \ts and \ts $i \in \Des(w)$, then \ts $c^w_{u,v}\ts = \ts 0$.

\smallskip
{\small $(1)$}  If \ts $i \notin \Des(v)$ \ts and \ts $i \notin \Des(w)$, then \ts $c^w_{u,v} \ts = \ts c^{ws_i}_{us_i \ts , \ts v}$.

\smallskip
{\small $(2)$}  If \ts $i \in \Des(v)$ \ts and \ts $i \in \Des(w)$, then \ts $c^w_{u,v} \ts = \ts c^{w}_{us_i \ts , \ts vs_i}$.

\smallskip
{\small $(3)$}  If \ts $i \in \Des(v)$ \ts and \ts $i \notin \Des(w)$, then
$$
c^w_{u,v} \, = \, c^{ws_i}_{us_i \ts , \ts v} \, + \, c^{w}_{us_i \ts , \ts vs_i}
\, + \, \ep(i,j,k) \. \sum_{(j,k)} \. c^{w}_{u t_{jk} \ts , \ts vs_i}.
$$
Here the summation is over all \ts $1\le j <k \le n$ \ts such that \ts $u(j)<u(k)$ \ts
and \ts $|\{j,k\}\cap \{i,i+1\}|=1$, and we set
$$
\ep(i,j,k) \, = \, \left\{\aligned
\. 1 & \quad \text{if} \ \ j=i \ \ \text{or} \ \  k=i+1 \\
\. -1 & \quad \text{if} \ \ k=i \ \ \text{or} \ \ j=i+1 \\
\. 0 & \quad \text{otherwise.}\endaligned\right.
$$
\end{lemma}

\smallskip

It is important to note that in all four cases of Knutson's
recurrence where \ts $(u,v,w) \to (u',v',w')$ \ts in the lemma,  we have
\ts $\inv(u') \ge \inv(u)+1$.  Thus, after iterating the recurrence,
it can stop only at \ts $u=\bw_\circ\ts$ or at zero terms.  In the former case we must
also have \ts $w =\bw_\circ$, and dimension equation~$(\oplus)$  gives \ts
$v=\one$.  Of course, we then have \ts $c^{\bw_\circ}_{\bw_\circ\ts,\ts \one} = 1$.

% \Colleen{Previous: ``Thus, after iterating the recurrence,
% it can stop only at \ts $u=\bw_\circ\ts$." But we can also stop at $0$!}
%
% Sure.  Good.

It is also worth noting that when \ts $\inv(u') > \inv(u)+1$, we have either
\ts $c^w_{u,v}=0$ \ts or \ts $c^{w'}_{u',v'}=0$ \ts by the equation~$(\oplus)$.
Therefore, \ts $\inv(u') = \inv(u)+1$ \ts is the only nontrivial possibility,
and the total number of steps in the iteration of the recurrence is exactly
\ts $\inv(\bw_\circ)-\inv(u) = \binom{n}{2} - \inv(u)$.  In other words, after
\ts $\binom{n}{2} - \inv(u)$ \ts iteration steps, we obtain a signed
summation of Schubert coefficients \ts $c^{\bw_\circ}_{\bw_\circ\ts,\ts \one}$,
as all other terms in the summation are equal to zero.

In summary, every Schubert coefficient \ts $c^w_{u,v}$ \ts with $u,v,w$ \ts
satisfying the dimension equation~$(\oplus)$ is equal to the number of positive
minus the number of negative terms in the summation obtained by iterating
Knutson's recurrence for \ts $\binom{n}{2} - \inv(u)$ \ts steps.  Our signed
puzzle rule is a reworking of a signed combinatorial interpretation
given implicitly by this algorithm.

%\newpage

\medskip

\section{The construction}\label{s:puzzle}

\subsection{The region} \label{ss:puzzle-region}
Consider a parallelogram shaped region $\Ga$ of size \ts $n\times \ell$, where
$\ell = \binom{n}{2}-\inv(u)$ as in the introduction.  We bicolor the region into
equilateral triangles as in  Figure~\ref{f:palette} below.

\begin{figure}[hbt]
	\includegraphics[width=7.8cm]{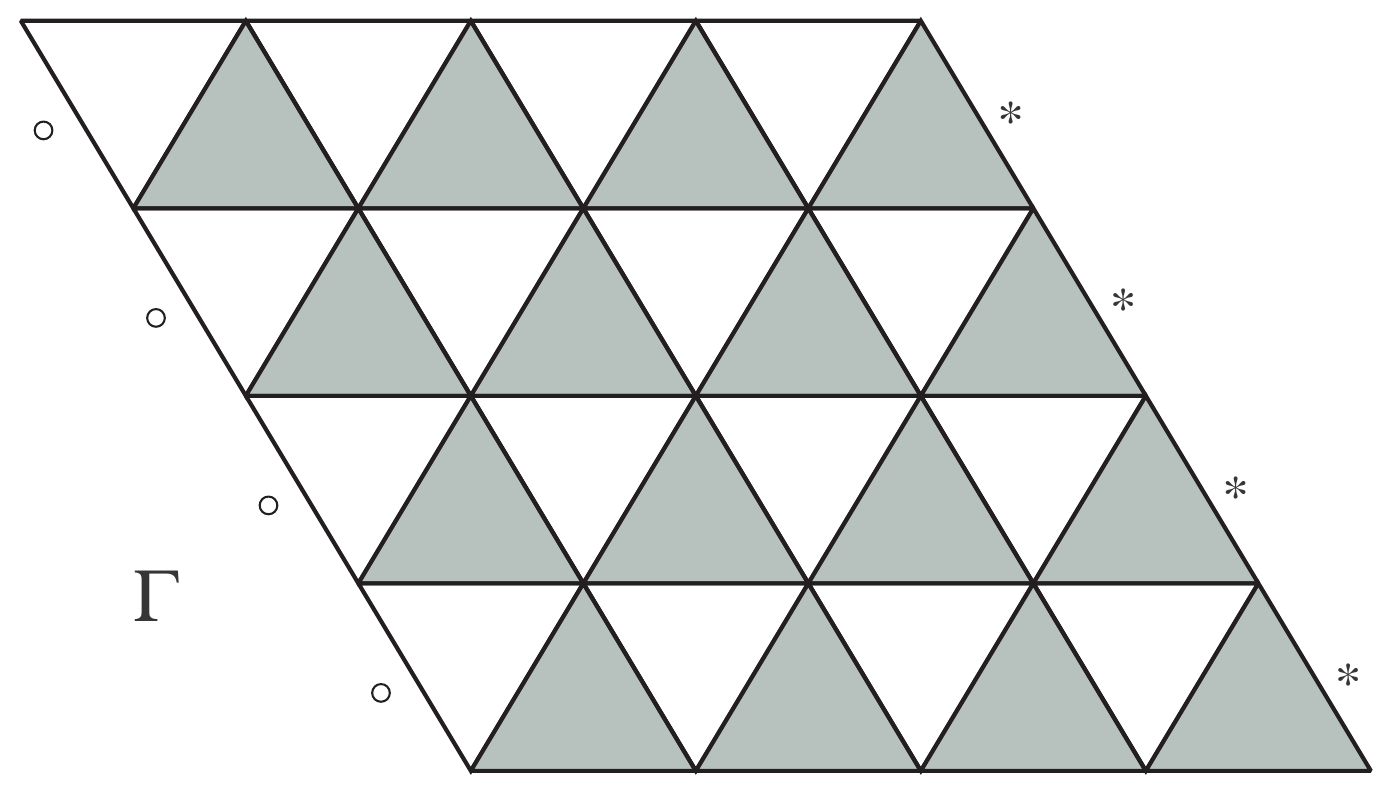}
%	\vskip-.25cm
	\caption{Region $\Ga$. }
	\label{f:palette}
\end{figure}

We label the boundary of $\Ga$  as follows.  Label all horizontal
edges with triples \ts $(a,b,c)\in [n]^3$.  The top edges are labeled
$(u(1),v(1),w(1))$, \. \ldots \. ,\. $(u(n),v(n),w(n))$.
The bottom edges are labeled
\ts $(n,1,n)$, \. $(n-1,2,n-1)$, \. \ldots \. \.,\. $(1,n,1)$.
Left and right edges on the boundary of~$\Ga$ are marked~$\offMark$ and $\onMark$, respectively. We use the term \emph{position} to mean a particular triangle
in~$\Ga$.
\smallskip

\subsection{Puzzle pieces} \label{ss:puzzle-pieces}
All puzzle pieces are equilateral triangles of three types: \emph{white},
\emph{shaded} and \emph{dark}.  We use \emph{triangle} to mean puzzle piece,
since all pieces will be triangular. In some cases we use the more specific terminology of \emph{triangle tile} to refer to a triangle, to avoid confusion.
% \Colleen{I vote to rename `triangle' to `piece'.}

\begin{figure}[hbt]
	\includegraphics[width=14.2cm]{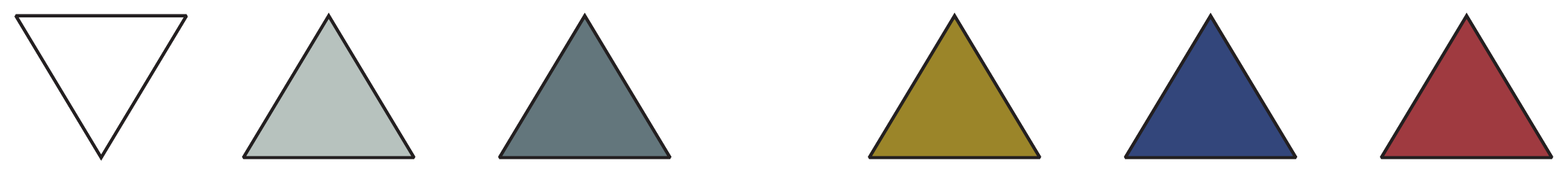}
	\vskip-.25cm
	 \caption{White, shaded and dark triangles.  Three colors of dark triangles: \\
dark yellow, dark blue and dark red.}
	\label{f:color}
\end{figure}

White triangle tiles are placed on white positions
in~$\Ga$, while shaded and dark triangle tiles are placed on shaded positions in~$\Ga$.  One should
think of dark triangle tiles as ``strongly shaded''; they encode a position
where Knutson's recursion is applied.  Dark triangles will come in three colors:
dark yellow, dark blue and dark red. The red and blue dark triangles correspond to positive and negative contributions,
respectively.  \emph{The sign of the puzzle will be the parity of the
number of red triangles in the puzzle}.
% Since the restrictions on red and blue triangles will be nearly identical,
To simplify definitions we will illustrate the dark triangles
as (uncolored) dark triangles, see Figure~\ref{f:color}.

% \Colleen{Igor: in grayscale, the dark yellow piece is too similar to the shaded piece. Could we lighten the gray in the 'shaded' piece or make the yellow darker? Ignore if too annoying.}
%
%  Done.

No rotations or reflections of the pieces are allowed, only parallel
translations.  The sides of the triangles will have labels and indicators,
described in Section~\ref{ss:puzzle-labels}.  Triangles are allowed to share a side in the puzzle if corresponding side labels and indicators are identical.

\begin{figure}[hbt]
	\includegraphics[width=10.2cm]{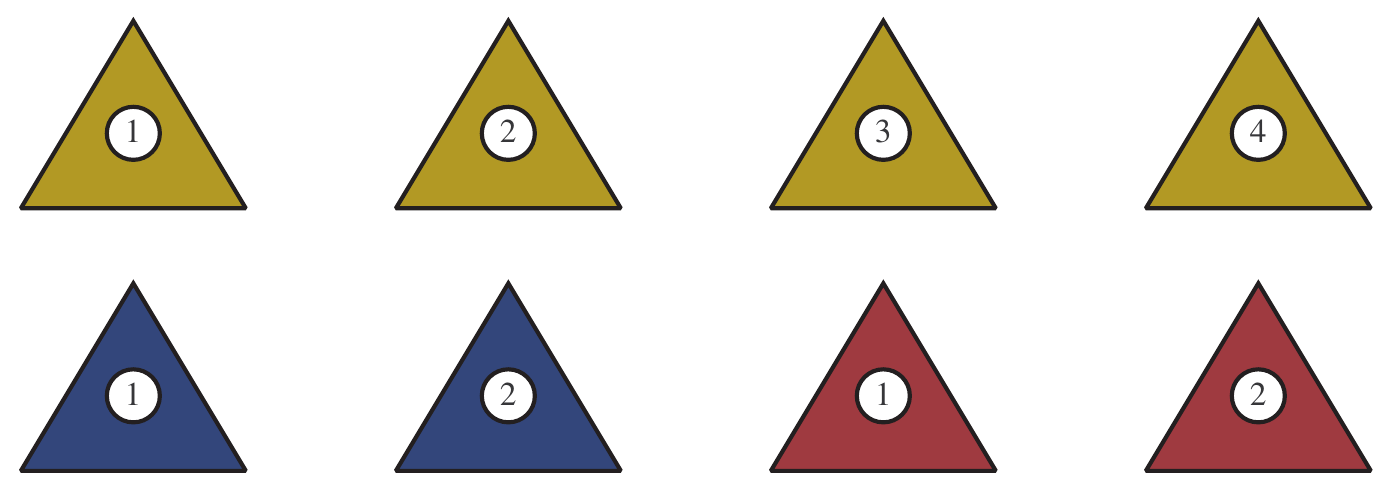}
	\vskip-.25cm
	\caption{Possible docket numbers of dark triangles. }
	\label{f:docket}
\end{figure}

Finally, in addition to color, all dark triangles have a \emph{docket number},
which ranges from~1 to~4 for yellow triangles, and is either~1 or~2 for
blue and red triangles, see Figure~\ref{f:docket}.  Docket numbers of yellow
triangles correspond to cases {\small (1)}, {\small (2)}, and first two
terms in {\small (3)} of Lemma~\ref{l:Knutson}.  Docket numbers of
blue and red triangles correspond to cases of positive and negative terms
in the summation in {\small (3)} of Lemma~\ref{l:Knutson}.
% These signs are given by $\ve(i,j,k)$.

\smallskip

\subsection{Labels and indicators} \label{ss:puzzle-labels}
Here, labels are numbers and indicators are symbols.  The labels and indicators on the triangles
will be somewhat involved and defined in stages.

% \Colleen{Igor: swap $\onMark$ with $\bullet$? Also small gripe, but the figure label looks like an O not $\circ$. Ignore if too annoying.}
%
%  Did the latter, not the former -- too much effort, and I like $\ast$.

\smallskip

\begin{figure}[hbt]
	\includegraphics[width=13.2cm]{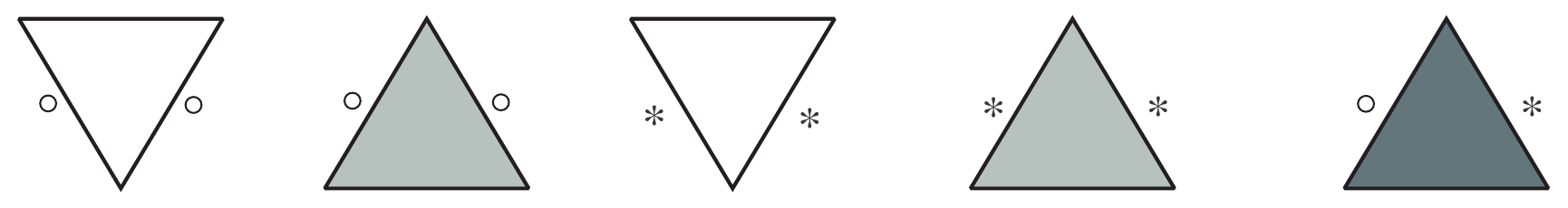}
	\vskip-.25cm
	\caption{Indicators~$\offMark$ and~$\onMark$ on the left and right edges of three types of triangles. }
	\label{f:LRmarking}
\end{figure}

\nin
\underline{Level~$0$.} \. We place an indicator $\offMark$ or $\onMark$ on the left and
right edges of all triangles as follows. For both white and shaded triangles, the indicators on the left and right edges must be equal.  For dark triangles, we mark the left edge with~$\offMark$
and right edge with~$\onMark$, see Figure~\ref{f:LRmarking}.  Since the left side of~$\Ga$ is
marked $\offMark$ and the right is marked~$\onMark$, these indicators ensure
that there is exactly one dark triangle in each row of~$\Ga$.

% \Colleen{change marking to `signal' or `indicator'?}
%done
\smallskip

\begin{figure}[hbt]
	\includegraphics[width=9.cm]{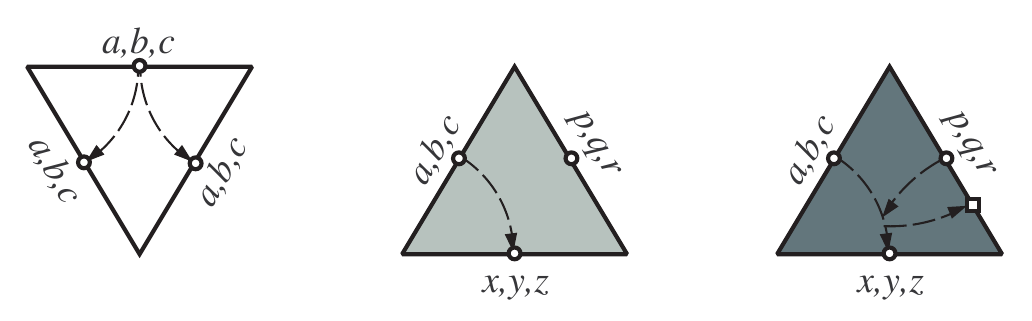}
	\vskip-.25cm
	\caption{Permutation labels on three types of triangles. }
	\label{f:triangle}
\end{figure}

\nin
\underline{Level~$1$.} \. All triangles have \emph{permutation labels}
on each edge.  These are triples $(a,b,c)$, where $a,b,c\in [n]$.
White triangles have three identical permutation labels. Shaded and dark triangles
can have distinct permutation labels on the left and right edges:  $a\ne p$, $b\ne q$,
and  $c\ne r$, see Figure~\ref{f:triangle}.  In the figures we often include
arrows to demonstrate the flow of permutation labels, i.e.\ how they shift through~$\Ga$.
These are not part of the labels and have only explanatory meaning.

\smallskip

\begin{figure}[hbt]
	\includegraphics[width=12.8cm]{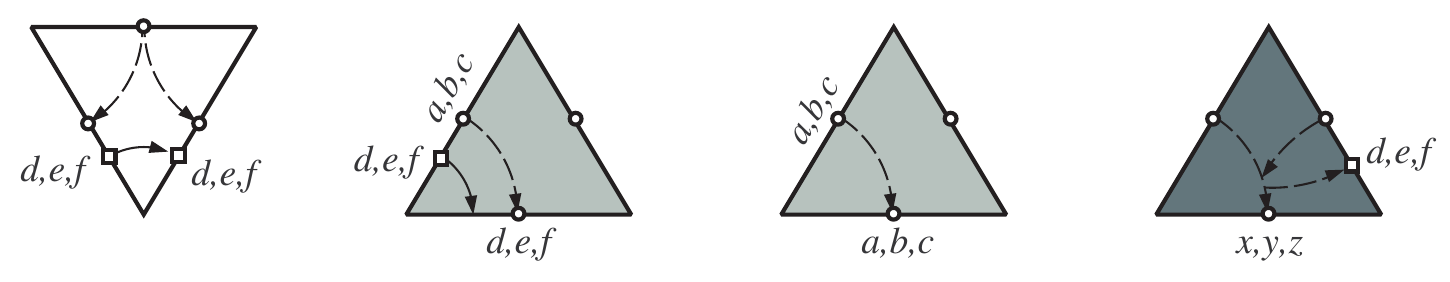}
	\vskip-.25cm
	\caption{Feedback labels on three types of triangles, and a shaded triangle
without feedback labels. }
	\label{f:feedback}
\end{figure}

\nin
\underline{Level~$2$.} \.  Some triangles have additional triples
of \emph{feedback labels} on the edges.  These will be of
the form $(d,e,f)$, where $d,e,f\in \<n\>$ and at least one of these is blank.
All dark triangles have feedback labels, which will appear only on their right edge.
Shaded triangles may have feedback labels, which will appear only on their left edge.
For shaded triangles \emph{with no} feedback nor transmuter labels
(see Level 3), the permutation labels on their left edge and bottom edge will be  equal.
Finally, white triangles may have equal feedback labels
on both left and right edges, see Figure~\ref{f:feedback}.
% \Colleen{Igor: please re-check Level 2 to ensure changes have introduced no errors}
\smallskip

\smallskip

\begin{figure}[hbt]
	\includegraphics[width=15.9cm]{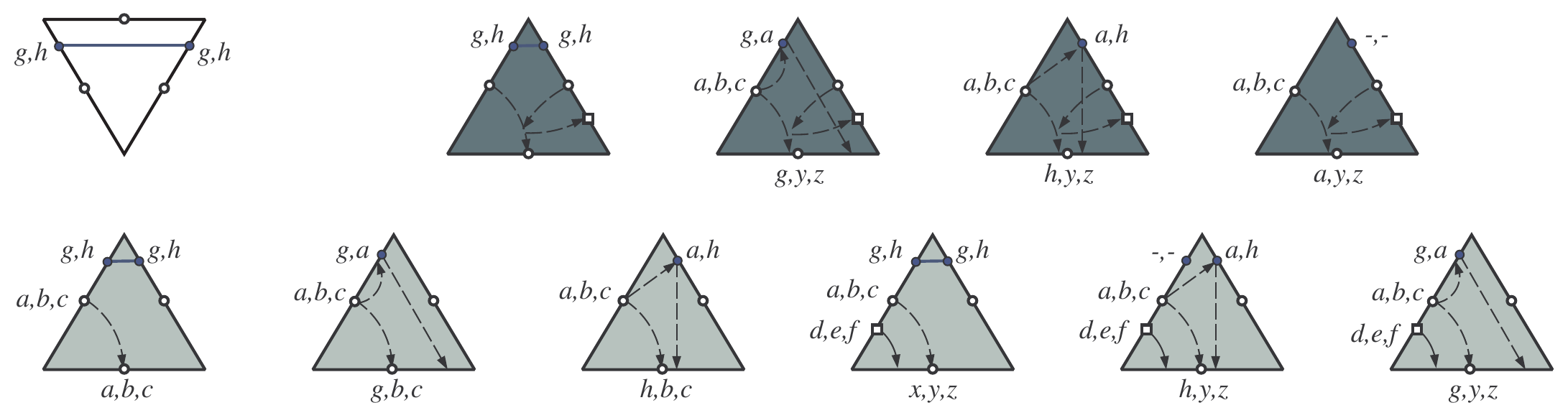}
	\vskip-.25cm
	\caption{Transmuter labels on three types of triangles. }
	\label{f:transmute}
\end{figure}

\nin
\underline{Level~$3$.} \. Finally, triangles may have additional
\emph{transmuter labels} on the left and right edges of the form $(g,h)$,
where $g,h\in \<n\>$.  In white triangles transmuter labels on the left and right edges must be equal.  Shaded and dark triangles
can have either equal transmuter labels on the left and
right edges, or have transmuter labels on only one edge. These transmuter labels can be
combined with permutation and feedback labels
as in Figure~\ref{f:transmute}.  For non-blank transmuter
labels $(g,h)$,  we always have the inequality $g<h$.

\smallskip

% Similarly, dark triangles have either equal
% transmuter labels on the left and right as in the figure.

\nin
{\it Note:} \ts For clarity, we distinguish the labels in the figures
by marking permutation labels with {\small $\circ$}\ts, the
feedback labels with \ts $\tinysquare\ts$, and the transmuter
labels with {\small $\bu\ts$}.  To see these markings, the reader
might want to zoom in.

% \Colleen{Igor: any meaning to different vertex styles? $\Box$ for feedback, $\circ$ for permutation, $\bullet$ for transmuter? This is not explicitly stated anywhere}

\medskip

\subsection{Constraints} \label{ss:puzzle-constraints}
We now describe constraints on the labels and indicators.
Roughly, there are very few additional constraints on white
and shaded triangles other than those described above.
Thus many types of labels can arise for white and shaded triangles.
On the other hand, the dark triangles are heavily constrained,
such that the docket number and color will uniquely determine
the constraints on the edges.

\smallskip

\begin{figure}[hbt]
	\includegraphics[width=14.8cm]{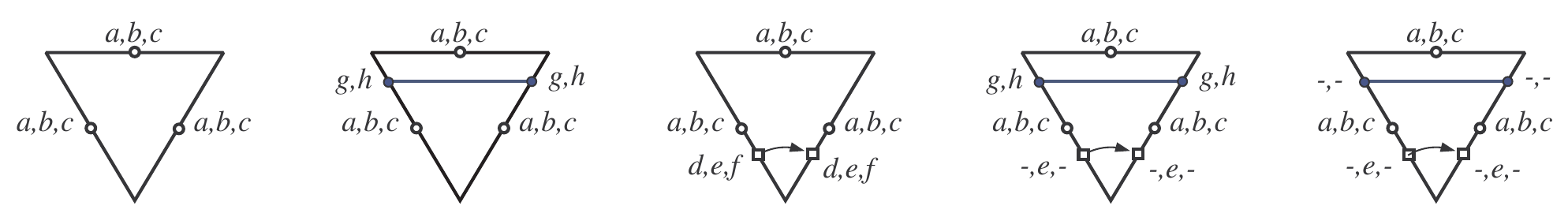}
	\vskip-.25cm
	\caption{Five types of white triangles. }
	\label{f:white}
\end{figure}

\nin
\emph{White triangles}:  There are five types of labelings of
white triangles depending on whether they have feedback labels,
transmuter labels, or both, see Figure~\ref{f:white}.
In the fourth and fifth triangles (from the left), two of the feedback
labels are blank.   Additionally, in the fifth triangle both transmuter
labels are blank.
For technical reasons, to align with the left boundary of~$\Ga$,
we allow the left permutation label to be empty.
There are two choices of the $\offMark/\onMark$ indicators of the first
two triangles, as described in Figure~\ref{f:LRmarking}.  We force that the last three triangles have~$\onMark$ on both edges.
This gives in total $O(n^6)$ white triangles.

\smallskip

\begin{figure}[hbt]
	\includegraphics[width=15.2cm]{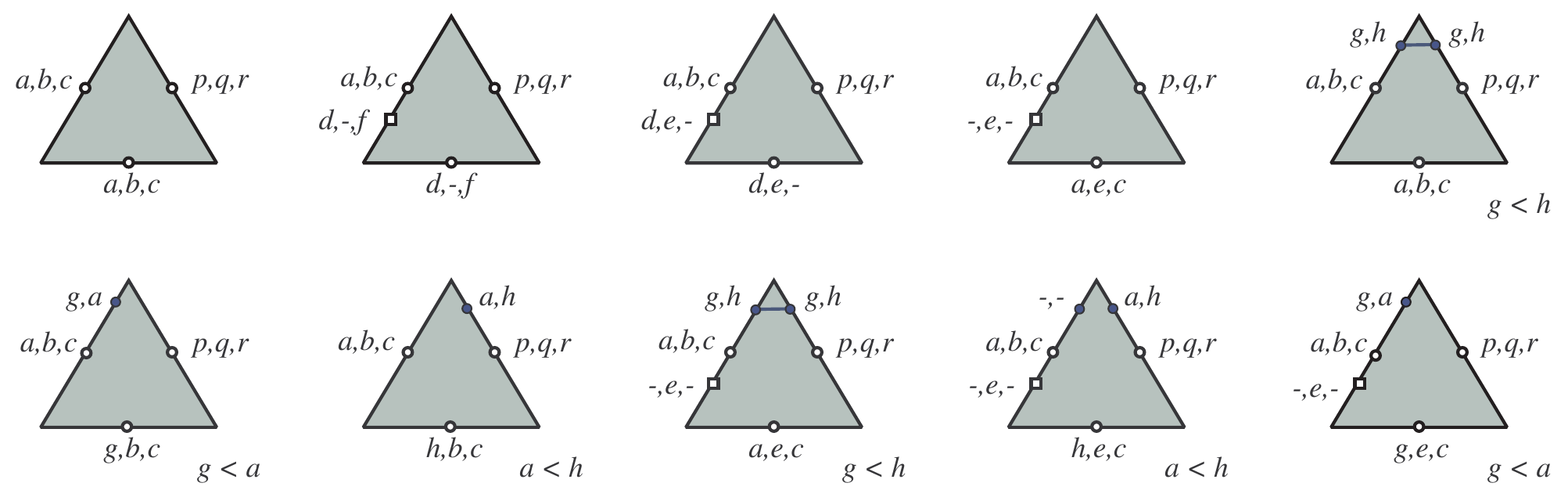}
	\vskip-.25cm
	\caption{Ten types of shaded triangles. }
	\label{f:shaded}
\end{figure}

\nin
\emph{Shaded triangles}: \ts There are ten types of labelings for
shaded triangles depending on whether they have feedback labels,
transmuter labels, or both, see Figure~\ref{f:shaded}.
For technical reasons, to align with the right boundary of~$\Ga$,
we allow the right permutation label to be empty.
% In each labeling there are two choices of $\offMark/\onMark$ indicators, as described in Figure~\ref{f:LRmarking}.
%
The third triangle in the second row has $9$ labels which
can all be distinct; other triangles have eight of fewer
distinct labels.  This gives in total $O(n^9)$ shaded triangles.

\smallskip

\begin{figure}[hbt]
	\includegraphics[width=9.6cm]{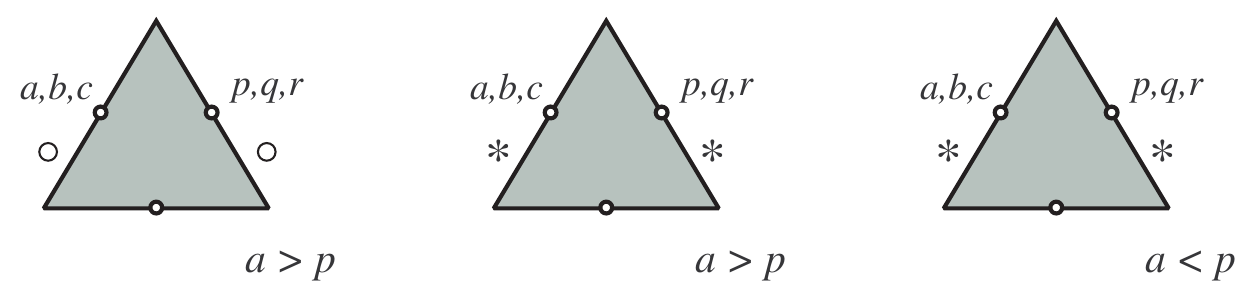}
	\vskip-.25cm
	\caption{Indicators of shaded triangles that are allowed under the ascent condition. }
	\label{f:star}
\end{figure}

% \nin
 For shaded triangles, there is an additional
condition relating indicators and permutation labels: a shaded triangle where the first entries in the permutation labels form an ascent
must have $\onMark$ indicators, see Figure~\ref{f:star}. We call this the \emph{ascent condition}. When these labels
form a descent, there is no additional constraint on the indicators.

\smallskip

\nin
\emph{Dark triangles}: \ts Recall all dark triangles have
indicators $\offMark$ on the left edge and $\onMark$ on the right.  Additionally,
they all have an ascent in the first entries of the permutation labels: $a<p$.  Below are
dark triangles arranged by color and docket number as in Figure~\ref{f:docket}.
Note that the additional inequalities for permutation labels correspond precisely to
the cases in Knutson's recursion.  This gives  $O(n^8)$ dark triangles in total.
% (of all docket numbers).

\smallskip

\begin{figure}[hbt]
	\includegraphics[width=12.5cm]{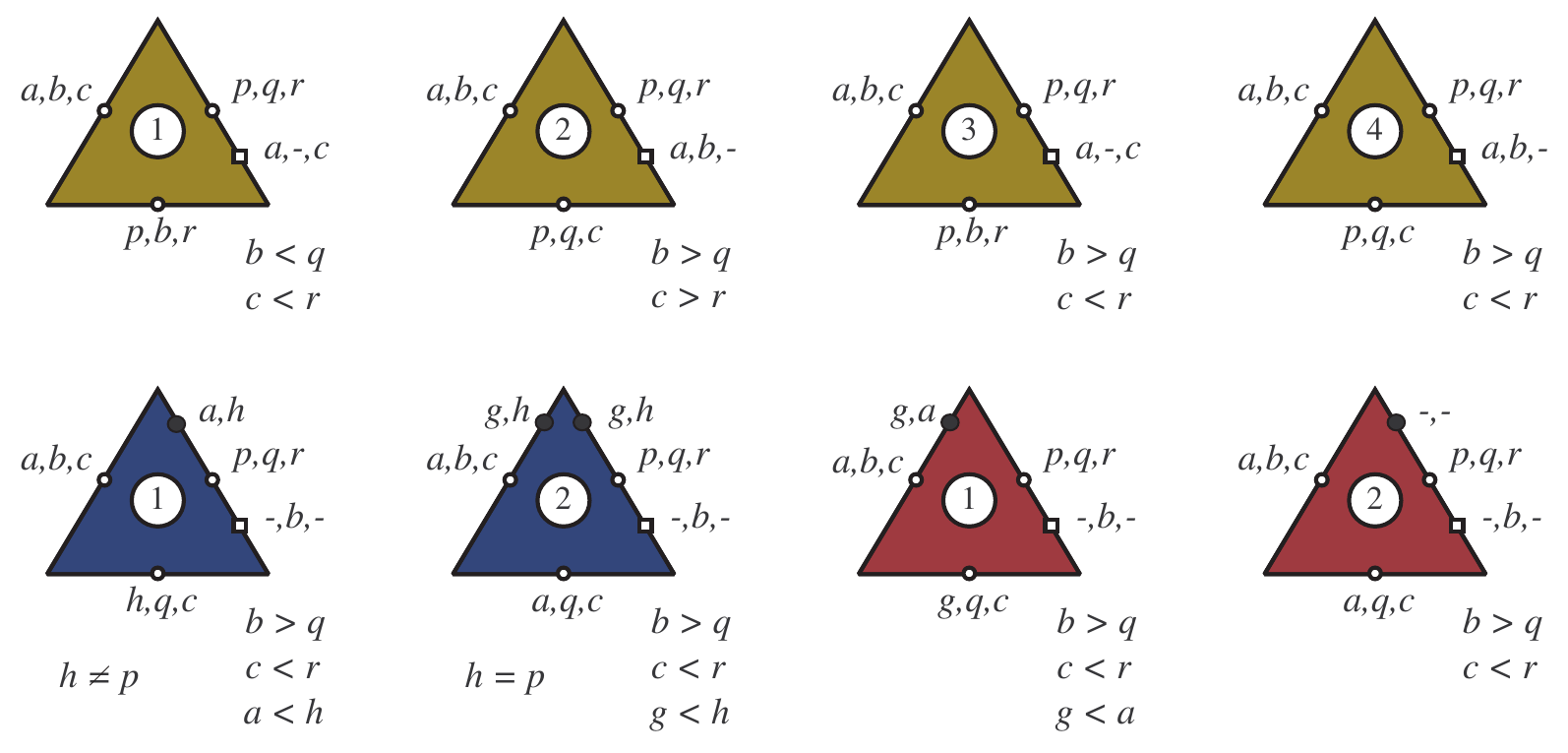}
	\vskip-.25cm
	\caption{Label constraints on dark triangles by docket number. }
	\label{f:dark}
\end{figure}

\smallskip

\subsection{Summary} \label{ss:puzzle-sum}
We gave a construction of $O(n^9)$ puzzle pieces.  All pieces are triangles
with three types of labels (permutation, feedback and transmuter),
two types of indicators ($\offMark$ and $\onMark$), five colors (white, shaded,
dark yellow, dark blue, and dark red), and docket numbers
to further distinguish dark colors.  We denote this set of triangle tiles
by $\cT_n$.

Given three permutations $u,v,w\in S_n$ which satisfy the dimension
equation $(\oplus)$, we constructed a parallelogram shaped region
$\Ga=\Ga(u,v,w)$ on a triangular grid with particular indicators and labels
on the boundary.  For each puzzle $T$ of~$\Ga$, the sign $s(T)$
is $(-1)^p$, where $p$ is the number of red triangles in $T$.
% \Colleen{n}
%
Denote by \ts $t_+(u,v,w)$ \ts and \ts $t_-(u,v,w)$ \ts
the number of puzzles $T$ with signs \ts $s(T)=1$ \ts and \ts $s(T)=-1$,
respectively.  Theorem~\ref{t:main} states that \ts
$c^w_{u,v} = t_+(u,v,w) - t_-(u,v,w)$.  We prove this in the next section.

\medskip

\section{Proof of Theorem~\ref{t:main}}\label{s:proof}

We prove the result by induction on the number $\ell$ of rows of~$\Ga$.
Thus, it suffices to show that the first row of triangles in the puzzle
is given by Knutson's recurrence as in Lemma~\ref{l:Knutson}.
Formally, we show that one step of the recurrence \ts $(u,v,w) \to (u',v',w')$ \ts
corresponds to top and bottom labels of triangles in the top row of a puzzle.

Recall there is exactly one dark triangle in each row of~$\Ga$,
and that the dark triangles correspond to index~$i$ in Lemma~\ref{l:Knutson}.
Then this dark triangle tile is placed the $i$-th shaded position
in its row of the region~$\Ga$.
By the constraints on dark triangles we must have \ts $i \notin \Des(u)$.
As mentioned before, the indicators $\offMark/\onMark$ on the boundary
of~$\Ga$ and on the triangles constrain the indicators which may appear.
% \textcolor{red}{there is a unique dark triangle
% with different indicators}.
Thus by the ascent condition on shaded triangle tiles, the index $i$ must be the
first ascent in the permutation~$u$.
% Note that labels of dark triangles are determined by their color and docket number.
% \Colleen{you assume the dark triangle corresponds to index $i$. So then do you need to suppose this is the $i$th triangle? Shouldn't this just follow from corresponding to index $i$? Do you mean the $i$-th shaded/dark triangle?}
% \Colleen{Igor: I've edited the above paragraph quite a bit. Please check}

We think of white triangles as splitters, which input the information,
i.e.\ permutation labels \ts $(u(i),v(i),w(i))$, \ts and transmit it to
triangles in shaded positions to its right and left.  When there are feedback or transmuter
labels, they transmit the signal with no changes from the left edge to the right
or vice versa, see Figure~\ref{f:white}.

\begin{figure}[hbt]
	\includegraphics[width=10.8cm]{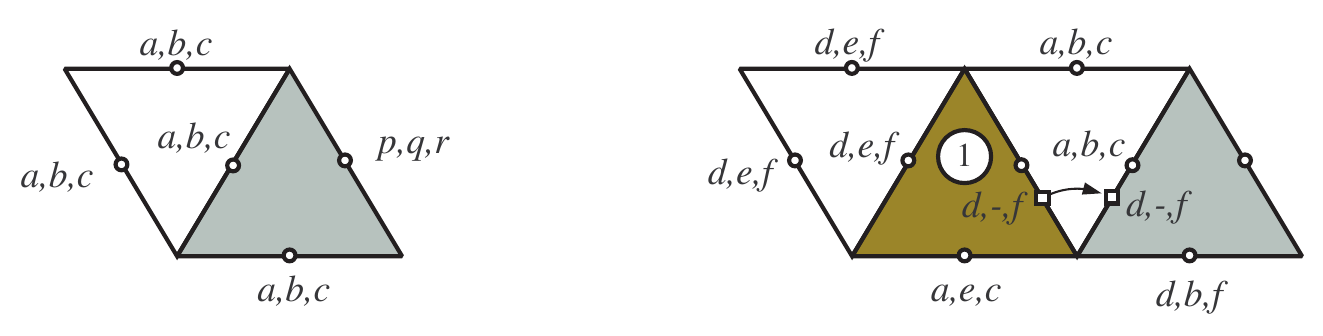}
	\vskip-.25cm
	\caption{Permutation labels transmitted from top to bottom. Feedback
signal dominated the permutation signal resulting in two transpositions:
$d\lra a$ and $f\lra c$.}
	\label{f:couple}
\end{figure}

Shaded triangles often play a similar role.  In particular, if no feedback or transmuter
labels appear, these triangles transmit the permutation labels from the left to the bottom edge.
In these cases the permutation labels \ts $(u(i),v(i),w(i))$ \ts are transmitted
unchanged from $i$-th top edge to $i$-th bottom edge of the first row of~$\Ga$,
see Figure~\ref{f:couple}.
% \Colleen{what about the other cases?}

\smallskip

The feedback labels in Figure~\ref{f:feedback} model transpositions
on the permutation labels. While the transpositions are initiated
in dark triangles, since shaded triangles are not adjacent to them,
the feedback labels transmit relevant permutation labels from left to right.
Thus when shaded triangles have feedback labels, these labels are on their left edge
to receive a signal sent by the dark triangle. This signal is transmitted by a neighboring white
triangle immediately to the left of the shaded triangle.  The feedback labels dominate
the permutation labels in determining the permutation label on the bottom edge of the shaded triangle; only blanks are substituted with permutation labels,
see Figure~\ref{f:couple}.

Dark triangles are distinguished by their color and docket numbers, each corresponding
to different cases of Knutson's recurrence.  The assumptions in these cases
are directly translated into constraints on the dark triangles given in
Figure~\ref{f:dark}.  For the dark yellow triangles, the transpositions
are local and the feedback labels enforce them.  However, for both the dark blue
and dark red triangles, this enactment of transpositions \ts $t_{jk}$ \ts can involve distant triangles.
Such transpositions are implemented with transmuter labels.

Transmuter labels are placed near the top of the triangle to be ``above the fray'',
see Figure~\ref{f:transmute}.
For some triangles (either dark or
shaded), they may affect permutation labels. However, in most cases, they introduce no constraints.  For transmuter labels $(g,h)$, we always
have label $g$ moving to the right, while $h$ to the left.
Similarly to feedback labels,  transmuter labels dominate the permutation labels in determining the permutation label on the bottom edge of the triangle.

By the dimension condition~$(\oplus)$, the recurrence {\small $(4)$} in
Lemma~\ref{l:Knutson} always involves transpositions which increase
the number of inversions in $u$ by one:
$$(\di) \qquad
\inv(ut_{jk}) \, = \, \inv(u) \. + \. 1.
$$
This gives a restriction $g<h$ for all transmuter labels.  This alone is not
a sufficient condition to ensure \ts $(\di)$\ts.
However, since the total number of rows in~$\Ga$ is \ts $\ell=\binom{n}{2}-\inv(u)$,
and since there is at least one inversion in \emph{every} row,
condition $(\di)$ holds automatically for all puzzles.

From this point, all conditions on dark blue and dark red triangles
are immediate translations of the last term of the summation {\small (3)}
in Lemma~\ref{l:Knutson}.   Transmuter labels can initiate at either dark or shaded
triangles.  The former possibility corresponds to having $j=i$ and dark
blue docket number~$1$, or having $k=i$ and  dark red docket number~$1$.
The latter corresponds to having $j=i+1$ and dark blue docket number~$2$,
or having $k=i+1$ and the dark red docket number~$2$.

\begin{figure}[hbt]
	\includegraphics[width=5.2cm]{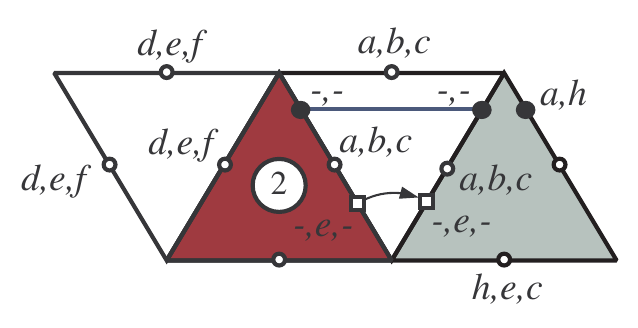}
	\vskip-.25cm
	\caption{For dark red triangle docket number~$2$, the blank transmuter labels signal initiates
transposition $t_{jk}$ with $j=i+1$ at the following shaded triangle. Here $d=u(i)$,
$a=u(i+1)$ and $h=u(k)$.}
	\label{f:blank}
\end{figure}

The only delicate point is the blank transmuter label when using
the dark red docket number~$2$ triangle.  This blank label signals that
the shaded triangle to the right should initiate the transpositions,
see Figure~\ref{f:blank}.
The details are straightforward.

\smallskip

Finally, we note that some constraints of the recurrence follow from the setup
and are never used.  Notably, we do not check that permutation labels $(u(1),\ldots,u(n))$,
$(v(1),\ldots,v(n))$ and $(w(1),\ldots,w(n))$, do indeed form permutations.  This is
assumed in the input, and for every row of horizontal edges this holds by induction.
Additionally, we never use the {\small $(0)$} case of the lemma.  This is because if we
have an ascent in both $u$ and $v$, it remains so by induction, and the desired puzzle
does not exist by the labeling of bottom edges of~$\Ga$.  This completes the proof. \qed
% \Colleen{explain more second to last sentence}
\smallskip

\begin{rem}\label{r:inv}
One can modify the construction to avoid the argument following equation
$(\di)$ above.  Note for $(\di)$ to hold, we must have $u(m)$ outside of the
interval \ts $[u(j),u(k)]$, for all \ts $j<m<k$.  This is a non-local constraint
which can be implemented by adding further constraints on the transmuter labels.
%
%increase in the number of inversions and the number of rows~$\ell$.
Formally, we can add the constraint that \ts $a\notin [g,h]$ on white triangles of the second
and fourth type in Figure~\ref{f:white}.
From above, this ensures that $(\di)$  holds automatically, without
referencing the number of rows~$\ell$ of~$\Ga$.  In the construction
above, we opted to avoid this modification for simplicity, but we will need
this version in the proof of Theorem~\ref{t:app} in Section~\ref{s:app}.
%  but at the cost of making
% \ts $\cT_n$ \ts  more cumbersome to define. % We omit the details.
\end{rem}
%\Colleen{Please reword (mainly first sentence). I can guess what you mean but this is hard to parse.}

\medskip

\section{Example}\label{s:ex}

Let $n=7$.  Take three permutations \ts $u=3251467$, \ts
$v=4126537$, \ts $w=6271534$ \ts in $S_7\ts$.  The following is an example
of Knutson's recursion steps with cases identified:

$$
\begin{bmatrix}
    3251647 \\
    4126537 \\
    6271534
  \end{bmatrix}  \xrightarrow{i=2, \ts (1)}
  \begin{bmatrix}
    3521647 \\
    4126537 \\
    6721534
  \end{bmatrix}  \xrightarrow{i=1, \ts (3b)}
\begin{bmatrix}
    5321647 \\
    1426537 \\
    6721534
  \end{bmatrix}  \xrightarrow[\ve(4,5,7)=-1]{i=4, \ts (3c)}
\begin{bmatrix}
    5321746 \\
    1425637 \\
    6721534
  \end{bmatrix}   \xrightarrow{i=4, \ts (1)}
$$
$$\begin{bmatrix}
    5327146 \\
    1425637 \\
    6725134
  \end{bmatrix}   \xrightarrow{i=3, \ts (1)}
  \begin{bmatrix}
    5372146 \\
    1425637 \\
    6752134
  \end{bmatrix}   \xrightarrow{i=2, \ts (2)}
    \begin{bmatrix}
    5732146 \\
    1245637 \\
    6752134
  \end{bmatrix}   \xrightarrow{i=1, \ts (1)}
  \begin{bmatrix}
    7532146 \\
    1245637 \\
    7652134
  \end{bmatrix}   \xrightarrow[\ve(5,4,7)=-1]{i=5, \ts (3c)} \
$$
$$
  \begin{bmatrix}
    7532164 \\
    1245367 \\
    7652134
  \end{bmatrix}   \xrightarrow{i=5, \ts (1)}
  \begin{bmatrix}
    7532614 \\
    1245367 \\
    7652314
  \end{bmatrix}  \xrightarrow[\ve(4,2,5)=1]{i=4, \ts (3c)}
  \begin{bmatrix}
    7632514 \\
    1243567 \\
    7652314
  \end{bmatrix}  \xrightarrow{i=4, \ts (1)}
 \begin{bmatrix}
    7635214 \\
    1243567 \\
    7653214
  \end{bmatrix}  \xrightarrow{i=3, \ts (2)} \quad
$$
$$ \begin{bmatrix}
    7653214 \\
    1234567 \\
    7653214
  \end{bmatrix}  \xrightarrow{i=6, \ts (1)}
  \begin{bmatrix}
    7653241 \\
    1234567 \\
    7653241
  \end{bmatrix}  \xrightarrow{i=5, \ts (1)}
  \begin{bmatrix}
    7653421 \\
    1234567 \\
    7653421
  \end{bmatrix}  \xrightarrow{i=4, \ts (1)}
  \begin{bmatrix}
    7654321 \\
    1234567 \\
    7654321
  \end{bmatrix} \, = \,
  \begin{bmatrix}
    \bw_\circ \\
    \one \\
    \bw_\circ
  \end{bmatrix}
   \hskip.45cm
$$

\vskip.4cm

\begin{figure}[hbt]
	\includegraphics[width=15.cm]{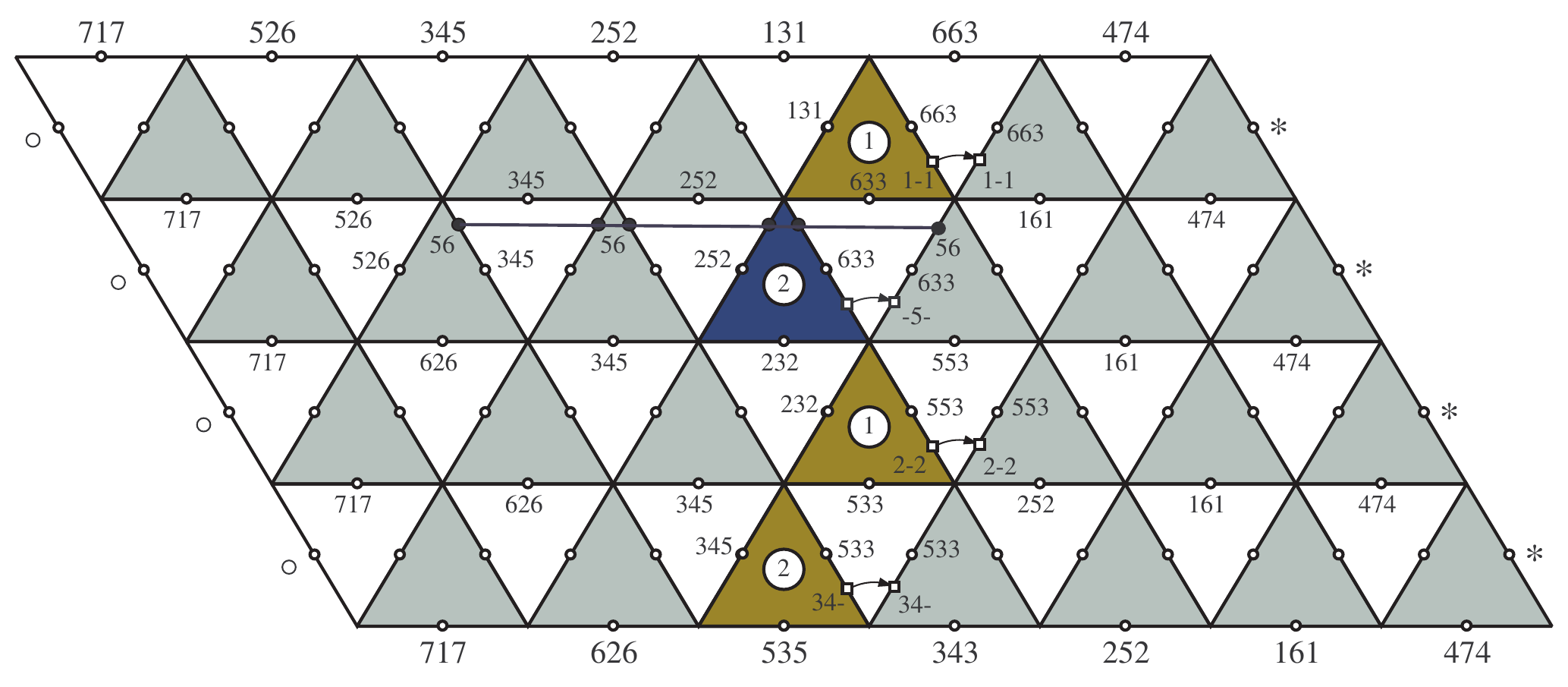}
	\vskip-.25cm
\caption{An example of a puzzle. }
	\label{f:example}
\end{figure}

\nin
Here $(3b)$ indicates the second summand in $(3)$ in Lemma~\ref{l:Knutson}, and so on.
The whole puzzle is quite large, so Figure~\ref{f:example} gives just the rows
corresponding to the third line of the calculation above, i.e.\ a quarter portion
of the actual puzzle.  To avoid cluttering we also omit some labels
which are clear from the example.

\medskip

\section{Proof of Theorem~\ref{t:app}}\label{s:app}
%
% \Colleen{for our conventions, we might want $\ov w = \bw_\circ\cdot  w\ts$}
% I think it's fine.  You might want to check again.
%

\subsection{The setup} \label{ss:app-setup}
First, note that \ts $c^w_{u,v} = c^{\ov u}_{\ov w,v}$, where
$\ov w = w\cdot \bw_\circ\ts$.  Thus we can rewrite
$$
\ga_k(n) \, = \, \sum_{u,v,w\ts\in\ts S_n \. : \. \inv(u)=\binom{n}{2}-k} \. c^w_{u,v}\..
$$
Now consider \emph{all} puzzles of the $n \times k$ parallelogram region $\Ga$ as
in Figure~\ref{f:palette}, where we remove the constraints on the top boundary of~$\Ga$.
Such puzzles contribute to \emph{some} triples of permutations $(u,v,w)$ as above.

Through the proof we will work with a modified set $\cT'_n$ of puzzle pieces given
in Remark~\ref{r:inv}.
Note that \ts $\cT'_n\subset \cT_n$ \ts by construction.
% with
% the difference given additional constraints on the labels that were not
% included in~$\cT_n$ for simplicity.
This is needed to
ensure \ts $\inv(u)=k$ \ts in every puzzle.  Then we have:
$$
\ga_k(n) \, = \, \sum_{\text{puzzle $T$ of $\Ga$ with $\cT_n$}} s(T)\ts.
$$

From this point on, to simplify the counting we will work with a rectangular
region obtained by an affine transformation of~$\Ga$ as in the figure below.
Here each position (of equilateral triangle shape) is turned into the right isosceles triangle,
so two such triangles (one white and one shaded) form a unit square.  We still refer
to the resulting region as~$\Ga$.  We now modify puzzle pieces, place them
on the new region accordingly, and refer to them in the same way as in the proof above.

\begin{figure}[hbt]
	\includegraphics[width=14.1cm]{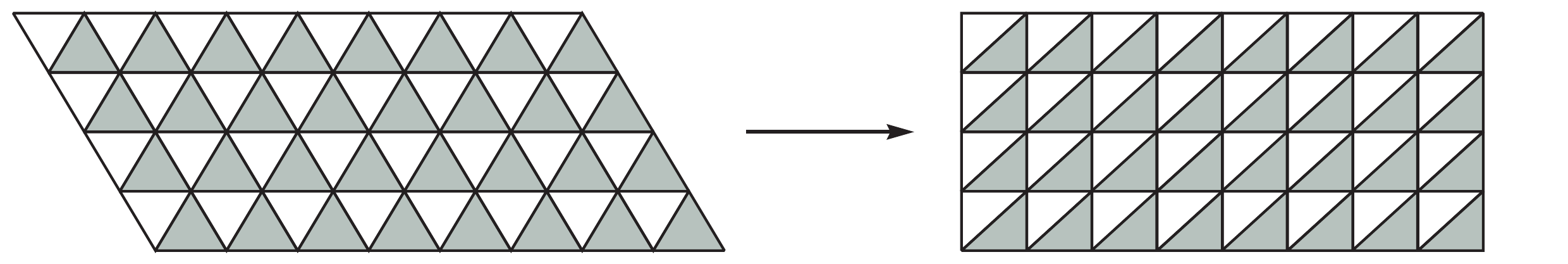}
	\vskip-.25cm
	\caption{Turning parallelogram region $\Ga$ into a rectangle.}
	\label{f:rectangle}
\end{figure}

\subsection{Relative placements and labelings} \label{ss:app-rel}
By construction, there is a finite number of types of dark triangles and
shaded triangles with a transmuter label on one side.
The positions and labelings of these triangles determines the puzzle.
Recall that there are exactly $k$ dark triangles and at most $2k$ such
shaded triangles, where at most~$k$ are not immediately following the dark triangles.
We call them \emph{separated shaded triangles}, or \emph{separated triangles} for short.

Since $k$ is fixed, the number of {relative placements} $\pi$ of these dark and
separated shaded triangles is also finite and depends only on~$k$ but not on~$n$.
Here by a \emph{relative placement} we mean
how these triangles are arranged in $\Ga$ relative to each other (above, to the left, to
the right, etc.), when one ignores the distances between them.  In other words, 
two puzzles $T,T^\ast$ have the same relative placement if they have the same number of 
dark and separated shaded triangles in each row, and $i$-th such triangle in $j$-th
row is to the left (right) of $i'$-th such triangle in $j'$-th row in~$T$ if and only
if the same holds in~$T^\ast$. 

\begin{figure}[hbt]
	\includegraphics[width=5.2cm]{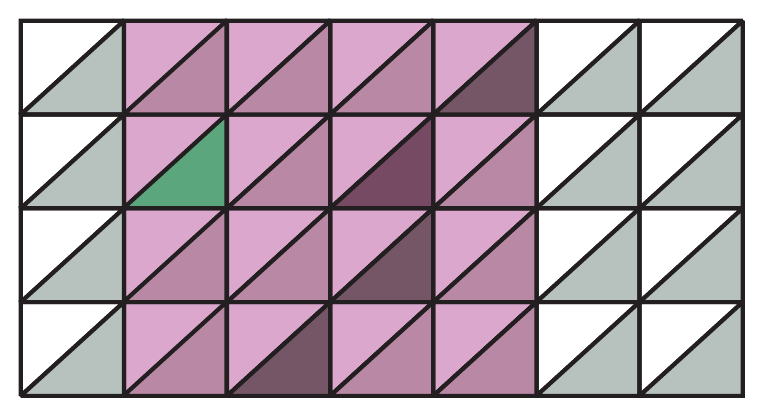}
	\vskip-.25cm
\caption{For the example in Figure~\ref{f:example}, we have \ts $I=\{2,3,4,5\}$,
and the columns are
highlighted purple.  Dark triangles are at $(1,5)$, $(2,4)$, $(3,4)$ and $(4,3)$.
Separated triangles at $(2,2)$ and $(2,5)$ are highlighted green. In this case \ts $I^\ast=\{1,2,3,4,5\}$. }
	\label{f:iset}
\end{figure}

Each relative placement
corresponds to \ts $f_\pi(n)$ \ts actual placements, where \ts $f_\pi$ \ts is
polynomial in~$n$. Here an \emph{actual placement} refers to a choice of positions
for the triangles that results in a puzzle.
One way to think of the actual placement is to think of the set  $I$ of columns
of~$\Ga$ which contains dark and separated triangles see Figure~\ref{f:iset}.
Clearly, $I$ determines the actual placement of these triangles, as their
relative positions to others determines their positions within the columns.

% \Colleen{It is slightly confusing that the separated triangle is yellow, when there are dark yellow triangles. Could we highlight purple or some other color that is not yellow, red, or blue?}
%
% Ugh...  Done.

We partition the set \ts $I \subseteq [n]$ \ts  into blocks
of consecutive integers: \ts $I = I_1 \cup I_2 \cup \ldots\ts$,
where \ts $I_1=[i_1,i_1']$, $I_2 = [i_2,i_2']$, etc.
Denote \ts $I^\ast: = I^\ast_1 \cup I_2^\ast \cup \ldots\ts$,
where \ts $I_m^\ast :=[i_m-1,i_m']$ \ts for all $m=1,2,\ldots$
Since \ts $|I|\le 2k$, we have \ts $|I^\ast|\le 4k$.
From this point, we will work with~$I^\ast$, see Figure~\ref{f:iset}.

Next, we need to take into account the number of possible
permutation, feedback and transmuter labels of \emph{all triangles} \ts which lie
in  columns~$I^\ast$.  The number~$r$ of such labels satisfies \ts
$r\le 9\cdot k \cdot |I^\ast| = O(k^2)$.
Given~$\pi$, there is a large number of equalities and
inequalities on these labels, resulting in polynomially many possible
labelings.  These labelings again can be characterized by the
\emph{relative labelings} determined by the at most \ts $r^2 = O(k^4)$ \ts
inequalities on the labelings.  We denote relative labelings by~$\la$.

\subsection{Actual labelings} \label{ss:app-count1}
Fix~$\pi$ and~$\la$ as above.  For all $n\ge 1$, 
denote by \ts $g_{\pi,\la}(n)$ \ts the number of possible actual
labelings with the relative placement~$\pi$ and relative labeling~$\la$.
As we show below, this is exactly the number of puzzles with given~$\la$ and~$\pi$. 

Although we will not need this fact, let us briefly show that
\ts $g_{\pi,\la}$ \ts is polynomial in~$n$.
Indeed,  \ts $g_{\pi,\la}(n)$ \ts is the number of integer
points in the $r$-cube $[n]^r$ minus some half-spaces of the type $x<y$,
restricted to hyperplanes of the type $x=z$, and outside of some hyperplanes
of the type $y=z$.
Resolving each $y\ne z$ as either $y>z$ or $y<z$, this shows
that \ts $g_{\pi,\la}$ \ts the sum of \emph{order polynomials} of a poset
on the labels, which is a polynomial in~$n$,
see e.g.\ \cite[$\S$3.15]{Sta99}.\footnote{More precisely,
some poset inequalities become strict; this does not affect the argument.}
% Let us emphasize that we use the fact that \ts $g_{\pi,\la}$ \ts depends only on the
% relative placement $\pi$, and not on the actual placement.
We will need a stronger argument of this type below.

From above, each such actual placement and actual labeling completely
determines a puzzle, except for the permutation labels not given
 by~$\la$.  Since the remaining permutation labels are unchanged
 from row to row by the dark and separated triangles, the bottom
 boundary determines them as well.  Indeed, these permutation labels
 in $i$-th column are given by \ts $u_i = w_i = n-i+1$ \ts and \ts
 $v_i = i$, for all $i \notin I^\ast$.

% In each block, the leftmost white triangles
% contain permutation labels from the left

Notice that in each block \ts $I_p^\ast=[i_p-1,i_p']$, the leftmost
white triangles in the column \ts $i_p-1$ \ts have equal permutation
labels in all rows: \ts $(n-i_p+2,i_p-1,n-i_p+2)$.  Similarly,
the rightmost shaded triangles in column \ts $i_p'$ \ts have
permutation labels on the right edges  \ts $(n-i_p',i_p'+1,n-i_p')$.
This implies that the actual labeling already contains the
information about~$I$, and thus actual placement of dark/separated
triangles.  In other words, relative placement $\pi$ and actual
labeling of triangles in columns $i \in I^\ast$ \ts uniquely
determines the whole puzzle.

%\Colleen{paragraph above is a bit hard to parse}

Of course, the inequalities on the actual labelings as above
can be inconsistent and the actual puzzle may not exist.
In summary, for every $\pi$ and actual labeling of triangles
contained in columns of dark/separated triangles,
there is  either one or zero possible puzzles.
Distinguishing between these possibilities is more difficult.

\subsection{Inequalities on actual labelings} \label{ss:app-ineq}
Fix~$\pi$ and~$\la$ as above.
We use parameters \ts $\al_1,\ldots,\al_r \in [n]$ \ts to denote
the actual permutation, feedback and transmuter labelings of all
triangles which lie in columns~$I^\ast$.
Whether the resulting puzzle exists or not introduces linear inequalities
on these parameters with integer coefficients as in Figures~\ref{f:shaded}
and~\ref{f:dark}, and the inequalities given by Remark~\ref{r:inv}.
%We need to be careful at this point.

First, as we mentioned above, the inequalities on the labelings coming
from each triangle are of the form \ts $\al_i < \al_j$, \ts $\al_i=\al_j$ \ts
and \ts $\al_i \ne \al_j\ts$. These give \ts $O(k^4)$ \ts inequalities.
Importantly, given $\pi$, $\la$ and $\{\al\}$, not all
permutation labels will form triples of permutations in each row.  Translating
this condition into permutation labels using relative orders $\pi$ and $\la$,
gives inequalities relating differences between the labels and distances
between the columns of their positions.  More precisely, we obtain
 inequalities of the form \ts $\al_i-\al_j > m$, \ts $\al_i-\al_j < m$ \ts
 or \ts  \ts $\al_i-\al_j = m$, which gives \ts $O(r^2) = O(k^4)$
 \ts additional inequalities.

Finally, we need to include the inequalities coming from constrains on the
transmuter labels, as  they relate to other labels given by Remark~\ref{r:inv}.
These inequalities are also of this type, but in logical combination.
Indeed, for example, for the shaded triangle that is
third in the second row on Figure~\ref{f:shaded}, the inequalities are
of the form the form \ts $a,p \notin [g,h]$.  This translates to
$$(\circledast) \qquad
\big((a<g) \vee (a>h)\big) \wedge \big((p<g) \vee (p>h)\big),
$$
where \ts $a,g,h,p \in \{\al_1,\ldots,\al_r\}$.
Note that the number of such inequalities on the parameters
is $O(k^2)$.  In total, the number of inequalities on the
labels is thus $O(k^4)$.

\subsection{Counting actual labelings} \label{ss:app-count2}
We now proceed to the counting of the set of labelings \ts $\{\al_i\}$ \ts
in \ts $[n]^r$ \ts which satisfy the inequalities as above.  The exact
inequalities or even their exact number will prove unimportant.
We will use only their form and the upper bound \ts $O(k^4)$ \ts of their number.

Fix a pair \ts $(\pi,\la)$ \ts of relative placements and labelings as before.
For simplicity, relabel all parameters $\al_i$ according to the relative
order~$\la$, so we have \ts $1\le \al_1 \le \ldots \le \al_r \le n$.
Denote by  \ts $J_{\pi,\la} \ssu [n]^d$
the set of possible vectors \ts $\bal=(\al_1,\ldots,\al_r)$ \ts as
above for which there exists a puzzle.

Resolve all inequalities $(\circledast)$
into two pairs of strict inequalities.  Similarly, resolve all inequalities
\ts $\al_i \ne \al_j$ \ts as either \ts $\al_i < \al_j$ \ts
or \ts $\al_i > \al_j$.  These define a partition of
\ts $J_{\pi,\la}$ \ts as a disjoint union of subsets $J_{\pi,\la,\vk}$,
where \ts $1\le \vk \le \ze(k)$ \ts and $\ze(k)=2^{O(k^4)}$.

Observe that \ts
$J_{\pi,\la} = n \ts Q_{\pi,\la} \cap \zz^r$, where the set \ts
$Q_{\pi,\la} \ssu [0,1]^r$ \ts is a disjoint union of \ts $\ze(k)$ \ts
convex polyhedra \ts $P_{\pi,\la,\vk} \ssu [0,1]^r$ with rational vertices.
Here each \ts $P_{\pi,\la,\vk}$ \ts is given by the integral inequalities
as above, where strict inequalities of the form  \ts $a<b$  \ts
are converted into nonstrict inequalities \ts $a\le b-1$.

Consider the \emph{Ehrhart quasi-polynomials}  \ts $h_{\pi,\la,\vk}(n):=|n \ts P_{\pi,\la,\vk} \cap \zz^r|$,
and observe that
$$
h_{\pi,\la}(n) \, := \, |J_{\pi,\la}| \, = \, |n \ts Q_{\pi,\la} \cap \zz^r| \, = \, \sum_{\vk=1}^{\ze(k)} \. h_{\pi,\la,\vk}(n)
$$
is also quasi-polynomial, see \cite[$\S$4.6]{Sta99} for the definitions.
Clearly, \ts $h_{\pi,\la}(n) \le f_\pi(n) \cdot g_{\pi,\la}(n)$.

% Note that increasing from $n$ to $(n+1)$ does not affect possibility of the puzzle.
% In other words, given \ts $(\pi,\la)$, either there is an actual placement and
% actual labeling for large enough~$n$, or such placement/labeling is impossible.

Now, write each polyhedron \ts $P=P_{\pi,\la,\vk}$ \ts by the defining inequalities as \.
$A \ts \bx \leqslant \bb$.  From above, every inequality can be rewritten in
the form \ts $\al_i - \al_j \le b$.  Thus, all maximal minors are determinants
of $r\times r$ \ts matrices with entries in \ts $\{0,\pm 1\}$, with
at most two $\pm 1$'s in every row of opposite sign.  Thus these minors are themselves in
\ts $\{0,\pm 1\}$ \ts by the same argument as in the standard proof of the
\emph{matrix-tree theorem}, see e.g.\ \cite[$\S$5.6]{Sta99}.
Therefore, all polyhedra \ts $P_{\pi,\la,\vk}$ \ts
are \emph{unimodular}, and thus have integral vertices, see e.g.\ \cite{Bar97}.
This implies that the quasi-polynomial \ts $h_{\pi,\la,\vk}(n)$ \ts is in fact a
polynomial in~$n$, for all~$(\pi,\la,\vk)$, ibid. %  and  the same holds for \ts $h_{\pi,\la}\ts$.

\subsection{Putting everything together} \label{ss:app-proof}
Observe also that the sign of a puzzle \ts $s(T)\in \{\pm 1\}$ \ts is
determined solely by $\pi$, so by a mild abuse of notation
we can write \ts $s(T) = s(\pi)$. Summing over all relative
placements and labelings, we have:
$$
\ga_k(n) \, = \, \sum_{(\pi,\la)} \.  s(\pi) \cdot  h_{\pi,\la}(n) \,
= \, \sum_{(\pi,\la)} \. \sum_{\vk=1}^{\ze(k)} \. s(\pi) \cdot  h_{\pi,\la,\vk}(n)\ts.
$$
The double summation has a constant number of terms for a fixed~$k$.  From above, each
\ts $h_{\pi,\la,\vk}(n)$ \ts is a polynomial in~$n$.  Thus, \ts
$\ga_k(n)$ \ts is also a polynomial in~$n$, as desired. \qed

\medskip

%\newpage
{\small

\section{Final remarks}\label{s:finrem}

\subsection{}  \label{ss:finrem-quote}
While discussing the background of puzzle rules in Schubert calculus,
Knutson and Zinn-Justin make the following observation:

\smallskip

\begin{center}\begin{minipage}{13.8cm}%
{{``We take $[$from above$]$ the oracular statement that
{\em puzzles should be related to Schubert calculus}.''}~\cite[p.~2]{KZ17}}\footnote{Original emphasis.}
\end{minipage}\end{center}

\smallskip

% \nin
% We  are somewhat less confident of puzzles' significance, as
% their omnipotence can be easily seen via equivalence with Wang tiles.

\nin
The signed puzzle rule in this paper is quite elaborate and uses a relatively
large number \ts $\Theta(n^9)$ \ts of puzzle pieces which are
not allowed to be rotated.  It is worth comparing this with some
of the earlier puzzle rules.

In the celebrated \emph{Knutson--Tao puzzles} \cite{KT03} for the
Littlewood--Richardson coefficients, there are only three puzzle pieces
and all $60^\circ$ rotations are allowed.  In the case of the
equivariant $K$-theory structure constants, Pechenik and Yong
\cite[Cor.~1.3]{PY17} modify and prove the previously conjectured
\emph{Knutson--Vakil puzzle rule}.  Their new puzzle pieces
can still be rotated, but now have complicated shapes
(this can be corrected by introducing new edge labels).
% Additionally, by their equivariant nature, their weights
% are polynomials and can be negative (in a controlled way).

%\Colleen{The signs in \cite{PY17} are indeed controlled (and inevitable) from K-theory. The formula is still cancellation-free, which should be noted. Citing \cite{KZ21} instead for an example of a cancellative puzzle formula makes more sense here.}

In the $3$-step case (for permutations with at most $3$ descents),
Knutson and Zinn-Justin gave several puzzle rules with the largest involving
$3591$ rhombi and some triangles, where now only $180^\circ$ rotations
are allowed \cite{KZ17,KZ21}.
In the $4$-step case, the number of puzzle pieces is even larger and some
of them have negative weight \cite{KZ21}.
Finally, in the \emph{separated descents} case, Knutson and Zinn-Justin
\cite{KZ23} have \ts $\Theta(n^2)$ \ts puzzle pieces.
%
%
% It is open as to whether such a puzzle rule exists.
% , although so did the construction by Knutson and Zinn-Justin in \cite{KZ23}.
We leave it to the reader to decide how our puzzles fit with these
earlier puzzle designs, and whether
this gives additional support to the quote above.
% We leave it to the reader to sort out philosophical
% implications, if there are any.

\subsection{}\label{ss:finrem-other}
There are several \emph{signed rules} for Schubert coefficients
known in the literature, sometimes in disguise. They are also called
\emph{signed combinatorial interpretations}, \emph{cancellative formulas}
and \emph{$\GapP$ formulas} in different contexts.  Perhaps, the cleanest signed
rule was given by Morales as a consequence of the Postnikov--Stanley
formula \cite[$\S$17]{PS09} and the pipe dream combinatorial
interpretation of Kostka--Schubert numbers, see \cite[$\S$10.2]{Pak-OPAC}.
Our own signed rule
in \cite[$\S$5]{PR24a} is somewhat similar but less explicit and stated
in a more general context.  More involved (and much more general)
signed rules are given by Duan \cite{Duan05} and Berenstein--Richmond
\cite{BR15b}.\footnote{The authors' insistence on using
{``Littlewood--Richardson coefficients''} to refer to Schubert coefficients
 is somewhat unfortunate as it initially obscures the very general nature
 of their results, see \cite[Remark~1.2]{BR15b}.}
Further generalizations of these rules are also known;
we refer to \cite{Knu22} for an overview.   While all these rules have
their own advantages, they seem incompatible with signed puzzle rules.

Additionally, there are several recursive formulas for computing
Schubert coefficients which can in principle be converted to formulas
for Schubert coefficients.  Notably, these include \emph{Billey's formula}
\cite[Eq.~(5.5)]{Bil99}, and the recent \emph{Goldin--Knutson formula} \cite{GK21}.
Unfortunately, Billey's formula requires a division, which is a major obstacle
to making it a signed rule.  In the case of the Goldin--Knutson formula,
the issue is the equivariant variables and their derivatives, which are
unavoidable, even in cases in which the dimension equation $(\oplus)$ holds.

% Neither of these are especially suitable for a signed puzzle rule.

\subsection{}\label{ss:finrem-inv}
It would be interesting to find a conceptual proof of Theorem~\ref{t:app}.
Can one compute the polynomials \ts $\ga_k$ \ts explicitly?  At the moment
we do not know even the degrees of \ts $\ga_k$ \ts beyond small special cases.
Curiously, the proof above only gives \ts $\deg \ga_k = O(k^2)$,
since the total number of labels of puzzle pieces can be rather large.
This is weaker than the elementary bound \ts $\deg \ga_k \le 6 k$ \ts
given in the introduction.

In a different direction, since Knutson's recurrence is originally stated
in the generality of equivariant cohomology, we expect that Theorem~\ref{t:app}
would also generalize in this direction.  It would be interesting to see if
the theorem generalizes to other cohomology theories mentioned in \cite{Knu22},
notably to $K$-theory and quantum cohomology.

\subsection{}\label{ss:finrem-ph}
To end on a philosophical note, it is worth pondering whether combinatorial
interpretations (rules), and, specifically, \emph{signed} \ts combinatorial
interpretations are worth studying.  As one would expect, there are
several schools of thought on the matter; see \cite[App.~B]{PR25} for
some background quotes.

In \cite[$\S$1.4]{Knu22},
Knutson lists three reasons why positive (unsigned) combinatorial
rules are better than signed: the vanishing problem,
computational efficiency and possibility of categorification.
We find the computational efficiency reason to be unconvincing,
at least from a theoretical point of view.  Indeed, conjecturally
the problem of computing Schubert coefficients is $\SP$-hard
\cite[Conj~1.2]{PR24b}.  Even if Schubert coefficients had a
$\SP$ formula, this formula might be quite hard to compute.
Since subtraction can be done in linear time, it is possible
and even likely that writing Schubert coefficients as the
difference of two $\SP$ functions can lead to a faster algorithm.
For example, famously, \emph{fast integer multiplication} \ts and
\emph{fast matrix multiplication} \ts algorithms
are fast \emph{because} they allow subtractions.  In other words,
if computing is the goal, then constraining oneself to positive
functions might not be a good strategy.

For the vanishing problem, we explain the state of art in our
recent papers \cite{PR24b,PR25}.
There, we obtain the best known vanishing results in full generality,
completely bypassing combinatorial arguments.
Unfortunately, even for the best studied $2$-step case, where a
puzzle rule was proved in \cite{BKPT16} (first conjectured by Knutson in~1999)
the complexity of the vanishing problem remains wide open.  Specifically,
in the $2$-step case, it would be very interesting to see if the
vanishing problem is in~$\P$ (this is known in the $1$-step case).
We are very far from resolving this problem despite having the puzzle rule.

In \cite{PR24a}, we adopt the opposite point of view, suggesting
that signed combinatorial rules have an intrinsic value, apart from
being a stepping stone towards a positive rule.  In fact, a
\emph{really good} \ts signed rule can be incredibly useful.
For example, the celebrated \emph{Murnaghan--Nakayama rule} for
the $S_n$ characters (see e.g.\ \cite[$\S$7.17]{Sta99}) is omnipresent
in the algebraic combinatorics literature, and has many powerful
applications in other areas.  This is despite the fact that its
absolute value has no positive rule, unless the polynomial
hierarchy collapses \cite{IPP24}.

This paper gives further evidence in favor of this point of view, as our
Theorem~\ref{t:app} gives a \emph{structural result} \ts that was
not easily attainable prior to the signed puzzle rule in Theorem~\ref{t:main}.
In fact, even the definition of \ts $\ga_k(n)$ \ts is mysterious from
algebraic point of view. However, fixing the height of the region
is completely natural in the tiling literature, see e.g.\ \cite{MSV00,Moo99}.
There, the number of tilings is usually computed using generating functions
and the \emph{transfer-matrix method}, see e.g.\ \cite[$\S$4.7]{Sta99}.
While the technical details are quite different, the connections
to Ehrhart theory have a similar flavor.

%\vfil

\vskip.6cm

% \small

\subsection*{Acknowledgements}
We are grateful to Sara Billey, Allen Knutson, Alejandro Morales,
Oliver Pechenik, Richard Stanley, Ada Stelzer, Anna Weigandt, Alex Yong and Paul Zinn-Justin
for interesting discussions and helpful comments.
This work was supported by
the National Science Foundation [CCF-2302173 to I.P., MSPRF DMS-2302279 to C.R.].
}

%\vfil
% \newpage

%\newpage
%\appendix

%\newpage

%%%%%%%%%%%%%%%%%%%%%%%%%%%%%%%%%%%%%%%%%%%%%%%%%%%%%%%%%%%%%%%%%%%%%%%%%%%%%%%%%%%%%%%%%%%%%
%%%%%%%%%%%%%%%%%%%%%%%%%%%%%%%%%%%%%%%%%%%%%%%%%%%%%%%%%%%%%%%%%%%%%%%%%%%%%%%%%%%%%%%%%%%%%
%%%%%%%%%%%%%%%%%%%%%%%%%%%%%%%%%%%%%%%%%%%%%%%%%%%%%%%%%%%%%%%%%%%%%%%%%%%%%%%%%%%%%%%%%%%%%
%%%%%%%%%%%%%%%%%%%%%%%%%%%%%%%%%%%%%%%%%%%%%%%%%%%%%%%%%%%%%%%%%%%%%%%%%%%%%%%%%%%%%%%%%%%%%
%%%%%%%%%%%%%%%%%%%%%%%%%%%%%%%%%%%%%%%%%%%%%%%%%%%%%%%%%%%%%%%%%%%%%%%%%%%%%%%%%%%%%%%%%%%%%
%%%%%%%%%%%%%%%%%%%%%%%%%%%%%%%%%%%%%%%%%%%%%%%%%%%%%%%%%%%%%%%%%%%%%%%%%%%%%%%%%%%%%%%%%%%%%
%%%%%%%%%%%%%%%%%%%%%%%%%%%%%%%%%%%%%%%%%%%%%%%%%%%%%%%%%%%%%%%%%%%%%%%%%%%%%%%%%%%%%%%%%%%%%
%%%%%%%%%%%%%%%%%%%%%%%%%%%%%%%%%%%%%%%%%%%%%%%%%%%%%%%%%%%%%%%%%%%%%%%%%%%%%%%%%%%%%%%%%%%%%
%%%%%%%%%%%%%%%%%%%%%%%%%%%%%%%%%%%%%%%%%%%%%%%%%%%%%%%%%%%%%%%%%%%%%%%%%%%%%%%%%%%%%%%%%%%%%
%%%%%%%%%%%%%%%%%%%%%%%%%%%%%%%%%%%%%%%%%%%%%%%%%%%%%%%%%%%%%%%%%%%%%%%%%%%%%%%%%%%%%%%%%%%%%

\end{document}